\documentclass[12pt]{amsart}
\usepackage{amsmath,amsfonts,amssymb,amscd,amsthm}
\usepackage[mathscr]{eucal}
\usepackage{relsize}
\usepackage{hyperref}
\DeclareMathOperator{\ind}{Ind}
\DeclareMathOperator{\restr}{Res}
\DeclareMathOperator{\sgn}{sgn}
\DeclareMathOperator{\tr}{tr}

\DeclareMathOperator{\sympsymbol}{Sp}

\DeclareMathOperator{\sorthsymbol}{SO}
\DeclareMathOperator{\unitsymbol}{U}

\DeclareMathOperator{\orthsymbol}{O}

\DeclareMathOperator{\id}{Id}

\newcommand{\fatPhi}{\mathbf{\Phi}}
\newcommand{\fatchi}{\pmb{\chi}}
\newcommand{\coloneqq}{:\hspace{0pt}=}

\newcommand{\reduc}[2]{\Big{\downarrow}^{#1}_{#2}}

\newcommand{\chiorthevenbd}[2]{\ensuremath{\fatchi^{\orthsymbol(2n)}_{#1}#2}}
\newcommand{\chiorthevenbdg}[1]{\chiorthevenbd{#1}{(g)}}
\newcommand{\chiortheven}[2]{\ensuremath{\chi^{\orthsymbol(2n)}_{#1}#2}}
\newcommand{\chisorthevenbd}[2]{\ensuremath{\fatchi^{\sorthsymbol(2n)}_{#1}#2}}

\newcommand{\chisortheven}[2]{\ensuremath{\chi^{\sorthsymbol(2n)}_{#1}#2}}

\newcommand{\chiorthoddbd}[2]{\ensuremath{\fatchi^{\orthsymbol(2n+1)}_{#1}#2}}

\newcommand{\chiorthodd}[2]{\ensuremath{\chi^{\orthsymbol(2n+1)}_{#1}#2}}
\newcommand{\chisorthoddbd}[2]{\ensuremath{\fatchi^{\sorthsymbol(2n+1)}_{#1}#2}}

\newcommand{\chisorthodd}[2]{\ensuremath{\chi^{\sorthsymbol(2n+1)}_{#1}#2}}

\newcommand{\chiorthbd}[2]{\ensuremath{\fatchi^{\orthsymbol(m)}_{#1}#2}}

\newcommand{\chisorthbd}[2]{\ensuremath{\fatchi^{\sorthsymbol(m)}_{#1}#2}}

\newcommand{\chisympbd}[2]{\ensuremath{\fatchi^{\sympsymbol(2n)}_{#1}#2}}
\newcommand{\chisympbdg}[1]{\chisympbd{#1}{(g)}}
\newcommand{\chisymp}[2]{\ensuremath{\chi^{\sympsymbol(2n)}_{#1}#2}}

\newcommand{\chiGbd}[2]{\ensuremath{\fatchi^{G}_{#1}#2}}

\newcommand{\schur}[2]{\ensuremath{{s}_{#1}{#2}}}
\newcommand{\schurbd}[2]{\ensuremath{\textbf{s}_{#1}{#2}}}

\newcommand{\myPbd}[2]{\ensuremath{\textbf{\textup{p}}_{#1} #2}}
\newcommand{\myPbdg}[1]{\myPbd{#1}{(g)}}
\newcommand{\myP}[2]{\ensuremath{p_{#1} #2}}

\newcommand{\sym}[1]{{\mathscr{S}_{#1}}}
\newcommand{\centr}[1]{{\mathscr{B}_{#1}}}

\newcommand{\shah}{Shahshahani}

\newcommand{\SBgroup}{{\sym{|\gamma|}\times\centr{k-|\gamma|}}}

\newcommand{\symp}[1]{\ensuremath{\sympsymbol ( #1 )}}
\newcommand{\sorth}[1]{\ensuremath{\sorthsymbol ( #1 )}}
\newcommand{\orth}[1]{\ensuremath{\orthsymbol ( #1 )}}

\newcommand{\unit}[1]{\ensuremath{\unitsymbol ( #1 )}}

\newcommand{\partition}{\vdash}

\newcommand{\ud}{\mathrm{d}}

\newtheorem{thm}{Theorem}
\newtheorem{lemma}{Lemma}
\newtheorem{prop}{Proposition}

\newcommand{\given}{\,\, | \,\,}
\newcommand{\esp}{\mathbb{E}}
\newcommand{\Complex}{\mathbb{C}}

\newcommand{\Integers}{\mathbb{Z}}

\begin{document}

\title[Averages over Lie groups and Weyl characters]{Averages over classical compact Lie groups and Weyl characters}
\author{Paul-Olivier Dehaye}

\thanks{This research was supported in part by the NSF grant FRG DMS-0354662.}

\date{\today}

\address{
 Department of Mathematics\\
 Stanford University\\
 CA
}
\email{pdehaye@math.stanford.edu}

\begin{abstract}
We compute $\esp_{G} (\prod_i \tr(g^{\lambda_i}))$, where $G=\symp{2n}$ or $\sorth{m}$ ($m=2n$, $2n+1$)  with Haar measure.  This was first obtained by Diaconis and \shah~\cite{DS}, but our proof is more self-contained and gives a combinatorial description for the answer. We also consider how averages of general symmetric functions $\esp_{G} \fatPhi_n$ are affected when we introduce a Weyl character $\chiGbd{\lambda}{}$ into the integrand. 
 We show that the value of $\esp_{G} \chiGbd{\lambda}{}\fatPhi_n / \esp_{G} \fatPhi_n$ approaches a constant for large $n$. More surprisingly, the ratio we obtain only changes with $\fatPhi_n$ and $\lambda$ and is independent of the Cartan type of~$G$. Even in the unitary case, Bump and Diaconis~\cite{BD} have obtained the same ratio. Finally, those ratios can be combined with asymptotics for $\esp_{G} \fatPhi_n$ due to Johansson~\cite{Jo} and provide asymptotics for $\esp_{G} \chiGbd{\lambda}{}\fatPhi_n$.
\end{abstract}


\maketitle

\section{Introduction}
Historically, the study of integrals of class functions over compact classical Lie groups with respect to Haar measure has been important for many areas of mathematics and physics. We will not even attempt to describe the relevance of this problem to physics, but refer the reader to the introduction of Mehta's book~\cite{Mehta}. On the mathematics side, we would like to mention at least the following works:
\begin{itemize}

\item The Heine-Szeg\"o identity and its relations to the strong Szeg\"o 
limit theorem. This identity expresses averages over unitary groups as 
determinants of Toeplitz matrices (see Bump and Diaconis \cite{BD}), while 
the strong 
Szeg\"o limit theorem gives asymptotics for such determinants (see the book by B\"ottcher and Silbermann~\cite{BS}).  

\item The study of averages of characteristic polynomials over compact classical Lie groups. Keating and Snaith conjectured that their calculations of those averages would serve as good predictors for moments of the Riemann $\zeta$ function \cite[unitary case]{KSzeta} and other data extracted from $L$-functions \cite[other classical groups]{KSLfunctions}. Our personal interest in Random Matrix theory sparks from this connection with Number Theory.
\item Diaconis and \shah's work \cite{DS} on averages of products of traces, and further refinements by Johansson \cite{Jo}. Those papers have a very probabilistic flavor, and rely on separate work for their most important result. Indeed, the answer to their computations turns out to be expressible as values of characters of the Brauer algebra. Those were evaluated by Ram~\cite{ArunRam,ArunRam2}, and are given by a rather complicated-looking function $g$ in \cite[Theorem 4]{DS}.
\end{itemize}

The first goal of this paper will be to offer with Theorem \ref{traces} a self-contained proof of the results of Diaconis and \shah, and even a combinatorial interpretation for the mysterious $g$ function that they obtained. If the reader only wants to understand the proof of this theorem, it might be helpful to observe that Propositions~\ref{derivation} and \ref{prime} include a $\gamma$ that will only be useful for Theorem~\ref{ratio}. The reader could thus safely assume that $\gamma=(0,0,\cdots)$ and still see a full proof of the following statement.
\begin{thm}
\label{traces}
Let $\lambda$ be a partition, $\lambda \partition k$ and $n \ge k$. 
Let $\epsilon = 1 $ when $G=\symp{2n}$ and $\epsilon=0$ when
$G=\sorth{2n}$ or $\sorth{2n+1}$. If 
$$
\myPbdg{\lambda} \coloneqq \prod_{i \in \mathbb{N}} \tr(g^{\lambda_i}) 
$$
then 
$$
\esp_{G} \myPbd{\lambda}{} = \sgn(\lambda)^\epsilon g(\lambda),
$$
where $g(\lambda)$ is defined to be the number of matchings of $k$ points
preserved under the action of a given element of $\sym{k}$ of cycle type
$\lambda$. 
\end{thm}
We remind the reader that a \emph{matching} of a set $S$ is a perfect partition of $S$ into pairs.

If we are willing to restrict the integrand to have $\lambda_i=1$ for all $i$, Rains~\cite[Theorem 3.4]{Rains} has proved this result in the full range for $n$. We present only the symplectic case of his result. In our notation, he proved that $\esp_{\symp{2n}} \myPbdg{\lambda}$ with $\lambda=(1,1,\cdots,1)\partition k$ is equal to the number of fixed-point-free involutions of length $k$ with no decreasing subsequence of length greater than $2n$. 
% ?? range

In the \emph{stable range}\footnote{See page \pageref{stable}.}, he is effectively counting the number of fixed-point-free involutions of length $k$, i.e. the number of matchings on $k$ points preserved by the identity permutation on those $k$ points.

The problem of Theorem~\ref{traces} was also solved in full generality by Pastur and Vasilchuk \cite{PV}, although their method of proof is arguably more complicated. We will sketch it in the orthogonal case. Let $F:\sorth{m}\rightarrow \mathbb{R}$ be a continuously differentiable function and $X$ be any $n\times n$ real antisymmetric matrix. By left-invariance of Haar measure, $\esp_{g \in \sorth{m}} F(e^{tX} g)$ is independent of the real parameter $t$ and so $\esp_{g \in \sorth{m}} (F'(g)Xg) = 0$, where $F'$ is the derivative of $F$. This expression can then be expanded and used to reduce the main expression to  simpler ones.

We would like to point out that our proof of Theorem~\ref{traces} involves the hyperoctahedral group $\centr{k}$. Both Stolz~\cite{Stolz} and Rains~\cite{RainsPhD} have already used the same group for this computation.

We now turn to a more complicated problem. 

Let $G$ be $\unit{n}$, $\sorth{2n}$, $\sorth{2n+1}$ or $\symp{2n}$ and let $\fatPhi_{n,f}$ be a class function on $G$, essentially defined by $\fatPhi_{n,f}(g)=\prod_i e^{f(t_i)}$, where $\{ t_i \}$ is a subset of eigenvalues of $g$ . There are extra technical conditions on $\fatPhi_{n,f}$, but these will be introduced just before the statement of Theorem~\ref{ratio}, Section~\ref{section3}. 

The strong Szeg\"o limit theorem gives the asymptotics and the rate of convergence of $\displaystyle \lim_{n \rightarrow \infty}(\esp_{\unit{n}} \fatPhi_{n,f})$. Johansson \cite{Jo} was the first to generalize this theorem to the other classical groups.

The second goal of this paper will be to study how those averages and asymptotics are affected when we introduce irreducible characters of $G$ into the integrand. The characters $\chiGbd{\lambda}{}$ were constructed by Weyl for the compact classical Lie groups using his Character Formula. By the Peter-Weyl theorem, these characters form a basis of the Hilbert space of class functions on $G$ and are thus very natural to consider.

Theorem~\ref{ratio} will show that the ratio $$\frac{\esp_G \chiGbd{\lambda}{}\fatPhi_{n,f}}{\esp_G \fatPhi_{n,f}}$$ approaches a limit when $n>>0$. This extends the corresponding results for the unitary groups due to Bump and Diaconis \cite{BD} to other classical groups. Remarkably, our ratio is independent of the Cartan type of the group $G$ and equal to the ratio they obtained for the unitary groups. It only varies with $f$ and $\lambda$ and can also be seen as the value achieved by the Schur polynomial $s_\lambda$ after setting the values of power polynomials to some Fourier coefficients of $f$.

A different point of view is offered in Bump, Diaconis and Keller \cite{BDK}: we can modify the Haar measure $\ud g$ into $\chiGbd{\lambda}{} \overline{\chiGbd{\lambda}{}} \ud g$. We know that $\chiGbd{\lambda}{} \overline{\chiGbd{\lambda}{}}$ is always positive and of mass 1 by orthogonality of Weyl characters hence $\chiGbd{\lambda}{} \overline{\chiGbd{\lambda}{}} \ud g$ is a measure. With this point of view, Theorem~\ref{ratio} would thus partially explain how the average of $\fatPhi_{n,f}$ with respect to Haar measure $\ud g$ is modified when \emph{twisting} the Haar measure by a character (see the last two remarks on page \pageref{remarks}).

% No justification of usefulness ???
% Clear it will be useful, but hard to say how
Thirdly, we would like to mention the recent preprint of Bump and Gamburd~\cite{BG}. They showed how many of the integrals useful for Number Theory can be computed in a unified way. An example of such an integral would be
$$
\int_{\unit{n}} \prod_i\Lambda_g(e^{\alpha_i})\, \ud g,
$$
where $\Lambda_g(\cdot)$ is the characteristic polynomial of $g$, and the 
$\alpha_i$'s are points on the unit circle. The importance of integrals of 
this type originates from the work of Keating and 
Snaith~\cite{KSLfunctions,KSzeta}, where the integrals have been shown to 
predict the moments of $\zeta(\cdot)$ and of L-functions.

The method of Bump and Gamburd is based on symmetric function theory and
classical results (Weyl Character Formula, Littlewood Branching Rules of
Theorem~\ref{branching}, page~\pageref{branching}, and Cauchy
Identity). The reader is referred to
their introduction for a much more comprehensive survey of all the results their method is known to produce, and how (if) they were proved before. 

This type of work is useful because it consolidates a wide array of methods into one more systematic technique. 

In that same vein, we hope that this paper can complement theirs to get closer to a more universal method. Indeed, we have shown how to introduce elements of the basis of symmetric functions into the integrand, an interesting step for that goal. Further steps are taken in the author's Ph.D. thesis and associated paper~\cite{DehayeBDTW}. 

Section~\ref{rep} will first go over notation, then introduce the reader to the representation theory of the compact classical Lie groups (Weyl characters and Branching Rules). Section~\ref{section3} will contain all of the proofs. It will also present the statement of Theorem~\ref{ratio}, and then shortly discuss its significance in relation to the rest of the literature. 

 The author is pleased to thank Daniel Bump and Persi Diaconis for numerous stimulating discussions. Alex Gamburd clarified some of the technical details of Section~\ref{Weyl} and suggested some of the references. 
% 
%The author would also like to acknowledge the people of the Prescott International Hotel and Hostel, Everett, MA, around whom this paper was drafted. 
%
 Finally, I would like to thank the people in my entourage for their unfaltering support.

\section{Representation theory of the classical groups}
\label{rep}
We now introduce Weyl characters and the branching rules between different 
classical compact Lie groups. We follow the expositions of \cite{BG} and 
\cite{KT}, but our notation is closer to \cite{BG}. 

\subsection{Notation}
\subsubsection*{Partitions. }
A partition $\lambda=(\lambda_1,\lambda_2,\cdots,\lambda_n)$ is a finite decreasing sequence of non-negative integers. We define the weight $|\lambda|$ of $\lambda$ to be the sum $\sum \lambda_i$. If this weight is $k$, we also use the notation $\lambda \partition k$. The length $l(\lambda)$ of $\lambda$ is the maximal $i$ such that $\lambda_i \neq 0$. The conjugate of $\lambda$ is denoted $\lambda'$. We say that a partition is even if all of its parts $\lambda_i$ are even. We define the union $\lambda \cup \mu$ to be the partition of $|\lambda|+|\mu|$ whose parts are the union of the parts of $\lambda$ and $\mu$. There is a partial ordering on partitions: $\lambda \subseteq \mu$ iff $\lambda_i\le \mu_i$ for all $i$. Finally, we define the $\lambda(i)$'s so that $(i^{\lambda(i)})=(\lambda_1,\lambda_2,\cdots,\lambda_n)$, i.e. $\lambda(i)$ counts the number of $\lambda_j$'s equal to $i$.

\subsubsection*{Symmetric group. }The symmetric group on $k$ points will be 
$\sym{k}$. If $\lambda \partition k$, elements of type $\lambda$ are the 
elements whose cycle types correspond to the partition $\lambda$. We use 
$\mathscr{C}_\lambda$ for the conjugacy class of those elements. We denote 
a centralizer in the group $G$ by $C_G(\cdot)$, and by $z_\lambda$ the 
order of the centralizer of an element of $\mathscr{C}_\lambda$. As usual, 
the irreducible characters $\chi_\lambda$ of $\sym{k}$ are indexed by 
partitions $\lambda \partition k$. We sometimes abuse notation and take 
$\chi_\lambda(\mu)$ to mean the value of $\chi_\lambda$ on 
$\mathscr{C}_\mu$. 
If $\chi_\lambda$ and $\chi_\mu$ are characters of $\sym{|\lambda|}$ and 
$\sym{|\mu|}$, their product $\chi_\lambda \odot \chi_\mu$ in the 
character ring of symmetric groups will be the character $\ind_{\sym{|\lambda|}\times 
\sym{|\mu|}}^{\sym{|\lambda|+|\mu|}}(\chi_\lambda \times \chi_\mu)$
(see Sagan's book \cite{Sagan} for all aspects of the representation theory of symmetric groups, and page 164 for the product of characters $\chi_\lambda \odot \chi_\mu$).

\subsubsection*{Classical groups. }Let $J$ be the $2n\times 2n$ matrix given by
$$
J=\begin{pmatrix}0 & -\id_n \\ \id_n & 0 \end{pmatrix}.
$$
We would like to introduce a few classical groups:
\begin{eqnarray*}
%\glinsymbol{(n,\Complex)} &=& \{ g \in M_n(\Complex) \given \det (g) \neq 0 \},\\
%\slinsymbol{(n,\Complex)} &=& \{ g \in M_n(\Complex) \given \det (g) = 1 \},\\
\unit{n} & = & \{ g \in M_n(\Complex) \quad \given \quad gg^*=I\},\\
\orth{n} &=& \{ g \in \unit{n} \quad \given \quad g g^t = I\},\\
\sorth{n} &=& \{ g \in \orth{n} \quad \given \quad \det(g) =1 \},\\
\symp{2n} &=& \{ g \in \unit{2n} \quad \given \quad gJg^t=J \}.\\
\end{eqnarray*}
If $G$ is one of those groups, it is compact for the topology induced by 
$M_n(\Complex)$ or $M_{2n}(\Complex)$. 
We can thus consider its Haar measure $\ud g$ and normalize it so the 
total volume of $G$ is 1. We write $\esp_G f$ for $\int_G f(g) \, \ud g$. 

\subsubsection*{Symmetric polynomials and power characters. } Let $\Complex[x_1,\cdots,x_m]^{\sym{m}}$ be the ring of symmetric polynomials in $m$ 
variables. We define the power polynomials 
$p_i(x_1,\cdots,x_m)=x_1^i+\cdots+x_m^i$ and 
$p_\lambda(x_1,\cdots,x_m)=\prod_i p_{\lambda_i}(x_1,\cdots,x_m)$. By abuse of notation, we also denote by $p_\lambda$ the generalized character of $\sym{|\lambda|}$ that is the indicator function with value $z_\lambda$ on the conjugacy class of type $\lambda$ (see Sagan~\cite[page 162]{Sagan}). The difference in the arguments of $\myP{\lambda}{}$ should prevent any ambiguity.
Note that the polynomial $p_\lambda$ is the image of the character  $p_{\lambda}$ under the 
characteristic map (see Bump's book \cite[Theorem~39.1]{Bump}). Finally, we 
define  the characters $\myPbd{\lambda}{}$ of $G=\unit{m},\orth{m},\sorth{m} $ or 
$\symp{m=2n}$ by $\myPbdg{\lambda}\coloneqq p_\lambda(t_1,t_2,\cdots,t_m)$ where 
the $t_i$'s are \emph{all} the eigenvalues of $g$. There is an obvious 
interpretation of those generalized characters in terms of the trace. For instance, we have $\myPbdg{(3,1,1)}=\tr (g^3)\cdot (\tr g)^2$.

\subsection{Weyl characters}
\label{Weyl}
Let $\lambda=(\lambda_1,\cdots,\lambda_n)$ be a partition. Let $i$ and $j$ 
be indices running between 1 and $n$. Guided by the Weyl Character 
Formula, we define the following symmetric functions of $\{x_1,\cdots,x_n\}$, 
actually polynomials in $\Integers[x_1,x_1^{-1},\cdots,x_n,x_n^{-1}]$:
\begin{eqnarray*}
\schur{\lambda}{(x_1,\cdots,x_n)} & = & \frac{\left|x_i^{\lambda_j+n-j}\right|}{\left|x_i^{n-j}\right|},\\
\chisorthodd{\lambda}{(x_1,\cdots,x_n)} &=& \frac{\left|x_i^{\lambda_j+n-j+1/2}-x_i^{-(\lambda_j+n-j+1/2)}\right|}{\left|x_i^{n-j+1/2}-x_i^{-(n-j+1/2)}\right|},\\
\chisymp{\lambda}{(x_1,\cdots,x_n)} &=& \frac{\left|x_i^{\lambda_j+n-j+1}-x_i^{-(\lambda_j+n-j+1)}\right|}{\left|x_i^{n-j+1}-x_i^{-(n-j+1)}\right|}.\\
\end{eqnarray*}
The $\schur{\lambda}{(\cdot)}$ are the regular Schur polynomials that 
appear in the representation theory of the symmetric group. Take $g \in 
\unit{n}$ (resp. $\sorth{2n+1}$ or $\symp{2n}$). Label the eigenvalues of 
$g$ by $\{ t_1,\cdots,t_n\}$ (resp. $\{t_1,t_1^{-1},\cdots , 
t_n,t_n^{-1},1\}$ or $\{t_1,t_1^{-1},\cdots , t_n,t_n^{-1}\}$). This 
allows us to define the functions $\schurbd{\lambda}{(g)}$, $\chisympbd{\lambda}{(g)}$ or $\chisorthoddbd{\lambda}{(g)}$ through the values of the respective function on the subset $\{t_1,\cdots,t_n\}$. 

When $G=\sorth{2n+1}$  (resp. $\symp{2n}$), Weyl showed that the character $\chisorthoddbd{\lambda}{}$ (resp. $\chisympbd{\lambda}{}$) is irreducible when $l(\lambda)\le n$. This is called the \emph{stable range} for $\lambda$\footnote{The book of Goodman and Wallach~\cite[Chapter~10]{GW} is the standard reference for this. See also the paper of Koike and Terada~\cite{KT}.}.
\label{stable}

Due to the involution in the Dynkin diagram of type $D_n$, the case of $\chisorthevenbd{\lambda}{}$ is actually special. We will again define $\chisorthevenbd{\lambda}{(g)}$ as the value of a function $\chisortheven{\lambda}{}$ on an appropriate subset of $n$ eigenvalues of $g$. The difference in this case is that we only have $\lambda_1\geq \lambda_2\geq\cdots\geq|\lambda_n|$ for the index set. If $\lambda$ is a regular partition, we define $\lambda_+ \coloneqq \lambda = (\lambda_1,\lambda_2,\cdots,\lambda_n)$ and $\lambda_- \coloneqq(\lambda_1,\lambda_2,\cdots,-\lambda_n)$. 
The characters $\chisorthevenbd{\lambda^+}{}$ and $\chisorthevenbd{\lambda^-}{}$ are exchanged by the involution on the Dynkin diagram, i.e by conjugation by an element of $\orth{2n}$ of negative determinant\footnote{It might be helpful for the reader to observe that in the odd orthogonal case, $\orth{2n+1}\cong \sorth{2n+1} \times \mathbb{Z}/2$ so the involution acts trivially.}. 

The Weyl character formula defines the functions
\begin{eqnarray*}
\begin{aligned}
\chisortheven{\lambda}{(x_1,\cdots,x_n)} =\hspace{2in}\\
\frac{\left|x_i^{\lambda_j+n-j}+x_i^{-(\lambda_j+n-j)}\right|+\left|x_i^{\lambda_j+n-j}-x_i^{-(\lambda_j+n-j)}\right|}{\left|x_i^{n-j}+x_i^{-(n-j)}\right|}.
\end{aligned}
\end{eqnarray*}
If we set $\chiortheven{\lambda}{}\coloneqq \chisortheven{\lambda^+}{} + \chisortheven{\lambda^-}{}$ when $\lambda_n\neq 0$ and $\chiortheven{\lambda}{}\coloneqq \chisortheven{\lambda^+}{}$ otherwise, then 
\begin{eqnarray*}
\chiortheven{\lambda}{(x_1,\cdots,x_n)} &=&
\frac{\left|x_i^{\lambda_j+n-j}+x_i^{-(\lambda_j+n-j)}\right|}{\left|x_i^{n-j}+x_i^{-(n-j)}\right|}.
\end{eqnarray*}
The character $\chiorthevenbdg{\lambda}$ is defined similarly by evaluating $\chiortheven{\lambda}{}$ on eigenvalues.
 
It is still a consequence of Weyl's work that $\chisorthevenbd{\lambda}{}$ 
is an irreducible character of $\sorth{2n}$ when $l(\lambda)\le n$. 
However, $\chiorthevenbd{\lambda}{}$ will merely be the character of the 
representation of $\sorth{2n}$ which is obtained by restricting an 
irreducible representation of $\orth{2n}$ to $\sorth{2n}$, \emph{not} the 
character of a representation of $\orth{2n}$. 

For the sake of uniformity in the orthogonal case, we will sometimes want to use 
$\chiorthoddbd{\lambda}\coloneqq\chisorthoddbd{\lambda}$. 

We also use the notational shortcut $\chiGbd{\lambda}$ where $G$ is one of the Lie 
groups defined above. It might be good at this point to remind the reader  that $\chi_\lambda$ denotes a character of a symmetric group.

The irreducibility of the various characters considered guarantees certain orthogonality properties, which we will only describe as needed in the proofs.

\subsection{Branching rules}
\label{branchingSection}
Let $G=\sorth{m}$ or $\symp{m}$. Since $G \subset \unit{m}$, the 
restriction of $\schurbd{\lambda}{}$ to $G$ is a class function for $G$ and can be expressed as a sum of $\fatchi^G_\mu$'s. The branching rules describe more precisely how to do that (see the paper of Koike and Terada \cite[page 492]{KT} for a modern and complete proof).

\begin{thm}[Littlewood]
\label{branching}
Let $\lambda$ be a partition of length less than or equal to $n$. Then 
\begin{eqnarray*}
\schurbd{\lambda}{}\reduc{\unit{2n}}{\symp{2n}} & = & \mathlarger{\sum_{\mu \subseteq \lambda}} \left(\sum_{\nu \textup{ even}} c^{\lambda}_{\nu'   \mu}\right) \chisympbd{\mu}{},\\
\schurbd{\lambda}{}\reduc{\unit{2n+1}}{\sorth{2n+1}} & = & \mathlarger{\sum_{\mu \subseteq \lambda}} \left(\sum_{\nu \textup{ even}} c^\lambda_{\nu  \mu}\right) \chiorthoddbd{\mu}{}\\
\schurbd{\lambda}{}\reduc{\unit{2n}}{\sorth{2n}} & = & \mathlarger{\sum_{\mu \subseteq \lambda}} \left(\sum_{\nu \textup{ even}} c^\lambda_{\nu \mu}\right) \chiorthevenbd{\mu}{},
\end{eqnarray*}
where $\schurbd{\lambda}{}\reduc{\unit{n}}{G}$ indicates the restriction to $G$ of the character $\schurbd{\lambda}{}$ of $\unit{n}$ and $c^\lambda_{\nu\mu}$ are the Littlewood-Richardson coefficients.
\end{thm}
\subsubsection*{Remark. } 
This is where the eigenvalue 1 "disappears" in the $\sorth{2n+1}$ case. Let $g\in\sorth{2n+1}\subset\unit{2n+1}$, with eigenvalues $\{1,t_1,$ $\cdots,t_n,$ $t_1^{-1},\cdots,t_n^{1}\}$. 
The left-hand side is $$\schurbd{\lambda}{(g)}=\schur{\lambda}{(1,t_1,\cdots,t_n,t_1^{-1},\cdots,t_n^{-1})},$$ while the right-hand side only involves terms of the form 
$$\chiorthoddbd{\mu}{(g)}=\chiorthodd{\mu}{(t_1,\cdots,t_n)}.$$

\section{Proofs}
\label{section3}
We will now present the main derivation. This is vaguely similar to a few steps of the proof of \cite[Theorem 2.1]{DE} in the unitary case. 
\begin{prop}
\label{derivation}
Let $\lambda \partition k$ and $n \ge k$ . Then 
$$
\esp_{\symp{2n}} \chisympbd{\gamma}{} \myPbd{\lambda}{} = \sum_{
\substack{\beta ' \textup{ even}\\\gamma \cup \beta \partition k}} \left<
\chi_{\gamma}{} \odot \chi_{\beta}{} , \myP{\lambda}{}
\right>_\sym{k}.
$$
Similarly (but with $\beta$ instead of $\beta'$), we have 
$$
\esp_{\sorth{2n+1}} \chisorthoddbd{\gamma}{} \myPbd{\lambda}{} = \sum_{
\substack{\beta \textup{ even} \\ \gamma \cup \beta \partition k}} \left<
\chi_{\gamma}{} \odot \chi_{\beta}{} , \myP{\lambda}{} \right>_\sym{k}
= \esp_{\sorth{2n}} \chisorthevenbd{\gamma}{} \myPbd{\lambda}{}
$$
Note: when $|\gamma| > |\lambda| = k$ or when $k-|\gamma|$ is odd, those sums
are indeed trivial and give a value of 0.
\end{prop}
\begin{proof}
The general method of proof is to use the branching rules from 
Section~\ref{branchingSection} to eventually transfer the problem to a
symmetric group. 

For definiteness, we will only prove this for $\symp{2n}$ and discuss at the end the minor changes needed in the orthogonal cases. Let $g \in 
\symp{2n}$ have eigenvalues $\{t_1,t_1^{-1},\cdots,t_n,t_n^{-1} \}$. Then
\begin{eqnarray*}
\myPbdg{\lambda}& =&\sum_{\mu \partition k}
\chi_\mu(\lambda) \schurbd{\mu}{(g)} \\
&=&\sum_{\mu \partition k} \chi_\mu{(\lambda)}\sum_{\nu \subseteq 
\mu}\left(\sum_{\beta ' \textup{ even}} 
c^\mu_{\nu
\beta} \right) \chisympbdg{\nu},\\
\end{eqnarray*}
where the first line follows from the usual decomposition of power polynomials 
into Schur polynomials given by the character table of a symmetric group (see Sagan~\cite[Equation (4.23)]{Sagan}). The second line follows by applying the branching rule for each $\mu\partition k$. The branching rule is only valid when $l(\mu) \le n$. This explains our final restriction of $n \ge k$.

We know that $\esp_{\symp{2n}} \chisympbd{\gamma}\chisympbd{\nu} = 1$ when
$\gamma = \nu$ and 0 otherwise (this is a consequence of the theory of the 
Weyl Character formula). Hence
$$
\esp_{\symp{2n}} \chisympbd{\gamma}{} \myPbd{\lambda}{} = \sum_{\mu
\partition k} \left(\chi_\mu(\lambda) \sum_{\beta' \textup{ even}}
c^\mu_{\gamma \beta}\right), 
$$
where the condition that $\nu = \gamma \subseteq \mu$ is still present
implicitly in the Littlewood-Richardson coefficient ( $c^\mu_{\gamma
\beta} =0 $ if $\gamma {\not \subseteq} \mu$). For the same reason, we see
that this sum is trivial when $|\gamma | > |\mu| = k$.

The final statement follows from observing that
$
\sum_{\mu \partition k} c^\mu_{\gamma \beta} \chi_\mu = \chi_\gamma \odot
\chi_\beta
$ 
and 
$
\chi(\lambda) = \left< \chi , \myP{\lambda}{} \right>_\sym{k}.
$

For the orthogonal groups, the only difference is that two characters will pop up when $\lambda_n \neq 0$. Let $m=2n$ or $2n+1$.
The branching rules will involve $\chiorthbd{\lambda}{}$ while the twist that we introduce comes from a character of type $\chisorthbd{\lambda}$.
Fortunately, all we need for the same proof to work is $\esp_G \chiorthbd{\lambda}{}\chisorthbd{\lambda}{}=1$:
\begin{eqnarray*}
\esp_G \chiorthevenbd{\lambda}{}\chisorthevenbd{\lambda}{}& = &\esp_G \chisorthevenbd{\lambda^+}{}\chisorthevenbd{\lambda}{} + \esp_G \chisorthevenbd{\lambda^-}{}\chisorthevenbd{\lambda}{} \\
&=&1+0 \text{ by orthonormality for $\sorth{2n}$.}\\
\esp_G \chiorthoddbd{\lambda}{}\chisorthoddbd{\lambda}{}& = &\esp_G \chisorthoddbd{\lambda}{}\chisorthoddbd{\lambda}{}\\
&=& 1 \text{ by orthonormality for $\sorth{2n+1}$.}
\end{eqnarray*}
\end{proof}

We would like to remind the reader at this point of a few facts from the
representation theory of the symmetric group. 
\begin{lemma}
Let $\sgn$ be the sign character in $\sym{k}$.
\begin{enumerate}
\item If $\beta \partition k$, then $\chi_{\beta'}{}$=$\sgn \otimes
\chi_{\beta}{} $,
\item If $\beta \partition k$, then $$\myP{\beta}{} \otimes \sgn =
\sgn(\beta) \myP{\beta}{}$$
\item Restrict $k$ to be even. Then 
$$\sum_{\substack{\beta \textup{ even}\\ \beta \partition k}}
\chi_{\beta}{} = \ind^\sym{k}_\centr{k}1, $$where $\centr{k}$ is the
centralizer of the chosen permutation $(1,2)$ $(3,4)\cdots(k-1,k)$ in
$\sym{k}$.
\item 
\label{useful4}
Restrict $k$ to be even. Then
$$
\sgn \otimes \ind_\centr{k}^\sym{k} 1 =
\ind_\centr{k}^\sym{k}(\restr_\centr{k}^\sym{k} \sgn).
$$
\end{enumerate}
\label{useful}
\end{lemma}
\begin{proof}
\begin{enumerate}
\item This is in Bump's book~\cite[Theorem 39.3]{Bump}. \label{conjugate}
\item This is immediate:
\begin{equation*}
\begin{split}
\myP{\lambda}{}\otimes \sgn & =  \sum_{\mu\partition k}
\chi_\mu(\lambda) \chi_{\mu}{}\otimes\sgn = \sum_{\mu\partition k}
\chi_\mu(\lambda) \chi_{\mu'}{}\\
\intertext{and, by part \ref{conjugate},}
&=\sum_{\mu\partition k} \chi_{\mu'}(\lambda) \chi_{\mu}{} 
=\sum_{\mu\partition k} \chi_\mu(\lambda)\sgn(\lambda)\chi_{\mu}
\end{split}
\end{equation*}
\item See \cite[Theorem 45.4]{Bump}.
\item This is a consequence of Frobenius Reciprocity.
% Consequence of inner product with a character 
\end{enumerate}
\end{proof}

This lemma leads immediately to a second version of Proposition~\ref{derivation}.
\begin{prop}
\label{prime}
Let $\lambda \partition k$ and $n \ge k$. 
Let $\epsilon = 1 $ when $G=\symp{2n}$ and $\epsilon=0$ when 
$G=\sorth{2n}$ or $\sorth{2n+1}$. 
Then
$$
\esp_G \chiGbd{\gamma}{} \myPbd{\lambda}{} =  \left<
\ind^\sym{k}_{\sym{|\gamma |} \times \centr{k-|\gamma |}}
\left(\chi_{\gamma}{}\otimes \sgn^\epsilon \right) ,
\myP{\lambda }{} \right>_\sym{k},
$$
where by a slight abuse of notation, we confuse $\sgn$ and
$\restr^\sym{k}_\centr{k} \sgn$.
\end{prop}
\begin{proof}
All the steps required are applications of Lemma \ref{useful} to the statement of Proposition \ref{derivation}. 
\begin{eqnarray*}
\esp_G \chiGbd{\gamma}{} \myPbd{\lambda}{} & = &
\sum_{
\substack{\beta \textup{ even}\\\gamma \cup \beta \partition k}} \left<
\chi_{\gamma}{} \odot (\sgn ^ \epsilon) \otimes \chi_{\beta}{} ,
\myP{\lambda}{} \right>_\sym{k} \\
&=&  \left< \chi_{\gamma}{} \odot \left( \sgn^\epsilon \otimes
\ind^\sym{k-|\gamma|}_\centr{k-|\gamma|} 1\right), \myP{\lambda}{}
\right>_\sym{k} 
\end{eqnarray*}
We now apply Lemma \ref{useful}.\ref{useful4} to get the result stated.

%\begin{eqnarray*}
%\phantom{\esp_G \chiGbd{\gamma}{} \myPbd{\lambda}{} }&=&  \left<
%\chi_{\gamma}{} \otimes \restr^\sym{k-|\gamma|}_\%centr{k-|\gamma|}
%\sgn^\epsilon, \restr^\sym{k}_{\sym{|\gamma|}\times\centr{k-|\gamma|}}
%\myP{\lambda}{} \right>_{\sym{|\gamma|}\times\centr{k-|\gamma|}} 
%\end{eqnarray*}
\end{proof}

\subsection{Discussion of Theorem~\ref{traces}}
As a special case to Proposition \ref{prime}, we are now ready to compute integrals of traces directly, without involving the Brauer algebra as in Ram~\cite{ArunRam2}. 
\begin{proof}[Proof of Theorem \ref{traces}]
We want here to compute $\esp_G \myPbd{\lambda}{}$, so we are now in the simplest case of Proposition $\ref{prime}$, when $|\gamma|=0$. 
When $k$ is odd, there is simply no matching on $k$ points. On the other 
hand, it was a
consequence of Proposition~\ref{derivation} that $\esp_{G}\myPbd{\lambda}{}=0$ as $k-|\gamma|=k$ is odd.
We can thus restrict our attention to the $k$ even case. We have thanks to Lemma \ref{useful} that
\begin{eqnarray*}
\esp_{G} \myPbd{\lambda}{}& = &  \left< \ind^\sym{k}_\centr{k} 1,
\myP{\lambda}{}\otimes \sgn^\epsilon \right>_\sym{k}\\
 & = & \sgn(\lambda)^\epsilon \left< 1, \restr^\sym{k}_\centr{k}
\myP{\lambda}{} \right>_\centr{k}\\
 & = & \frac{z_\lambda \sgn(\lambda)^\epsilon}{|\centr{k}|} \#
\left\{\sigma \in C_\sym{k}((1,2)\cdots(k-1,k)) \given \text{ type}(\sigma) = \lambda \right\}, 
\end{eqnarray*}
since $p_\lambda$ is an indicator function for the conjugacy class of
permutations of type $\lambda$ in $\sym{k}$.

If $\sigma \in C_\sym{k}((1,2)\cdots(k-1,k))$ then $\sigma$ preserves the matching $\{\{1,2\},$ $ \cdots , \{k-1,k\}\}$, i.e. it sends a pair to a pair. We use this to switch to the language of matchings.

\begin{eqnarray*}
\esp_{G} \myPbd{\lambda}{}& = &
\frac{\sgn(\lambda)^\epsilon}{|\mathscr{C}_{\lambda}|}
\frac{|\sym{k}|}{|\centr{k}|} {\# \left\{\sigma \in
C_\sym{k}((1,2)(3,4)\cdots(k-1,k)) \cap \mathscr{C}_\lambda \right\}}\\
&=&\frac{\sgn(\lambda)^\epsilon}{|\mathscr{C}_{\lambda}|}
\sum_{\text{matching $M$ of $k$ points}} {\# \left\{\sigma \in
\mathscr{C}_\lambda \given \sigma(M)=M\right\}}\\
&=&\frac{\sgn(\lambda)^\epsilon}{|\mathscr{C}_{\lambda}|} \#\{ (M,\sigma) \given M \text{ a matching of $k$ points},\\
& &\hspace{2in} \sigma \in \mathscr{C}_\lambda, \sigma(M)=M \}
\\
&=&\frac{\sgn(\lambda)^\epsilon}{|\mathscr{C}_{\lambda}|} \sum_{\sigma
\in \mathscr{C}_\lambda} \# \left\{ \text{matchings preserved by }  \sigma
\right\}.
\end{eqnarray*}
The last steps make use of a double-counting argument. All the summands in the last line are equal, and there are $|\mathscr{C}_\lambda|$ of them so we have
$$\esp_G \myPbd{\lambda}{}=\sgn(\lambda)^\epsilon g(\lambda),$$
where $g(\lambda)$ is the number of matchings preserved by a permutation of cycle type~$\lambda$.
\end{proof}

As mentioned earlier, this offers a combinatorial interpretation for a result first proved by Diaconis and \shah~\cite{DS}.
Naturally, we have to check that our definition of $g$ agrees with the definition they gave. This is a purely combinatorial problem.
\begin{prop}
Let $\lambda \partition k$. Then $g(\lambda)= \prod_j g_j(\lambda(j))$,
where $g_j(\cdot)$ is given by
\begin{eqnarray*}
\text{if $j$ is odd}\quad g_j(a) & = &
\left\{\begin{tabular}{ll}
0  &  \text{if $a$ is odd}\\
$j^{a/2} (a-1)(a-3)\cdots 1 $ & \text{if $a$ is even,}
\end{tabular}\right.\\
\text{if $j$ is even}\quad g_j(a)& = &\sum_t \binom{a}{2t}
 j^t (2t-1)(2t-3)\cdots 1  
\end{eqnarray*}
\end{prop}
\begin{proof}
From our combinatorial definition of $g$, it is immediate that
$g(\lambda) = \prod_j g((j^{\lambda(j)}))$. All we have left to prove is
$g((j^a)) = g_j(a)$.
\begin{description}
\item[if $j$ is odd: ] Take $\sigma \in \mathscr{C}_{(j^a)}$. Since each
cycle of $\sigma$ is of odd length, any matching of points preserved by
$\sigma$ must match cycles as well. If $a$ is odd there is no such
matching. If $a$ is even, any matching of points will also match
cycles. There are $(a-1)(a-3)\cdots 1$ possible matchings of cycles. Once a
matching of cycles is chosen, we still have to decide on how to match points
in each individual pair of cycles. There are $j$ choices for each of the
$a/2$ pairs of cycles. 
\item[if $j$ is even: ] This is more subtle, as matchings of points inside
the same cycle are allowed. Say there are $2t$ cycles whose points are
matched with points in another cycle (the \emph{external} cycles) and thus
$a-2t$ cycles whose points are matched with a point within the same cycle
(the \emph{internal} cycles). There are $\binom{a}{2t} $ ways
of choosing which cycles will be external, and then $(2t-1)(2t-3)\cdots 1$
ways of matching external cycles. Once we have a pair of external cycles,
there are $j$ ways of matching points between the two cycles. On the other
hand, there is a unique way of matching points within an internal cycle: a
point has to be paired with the point most distant for the ordering given
by the cycle.
\end{description}
\end{proof}

\subsection{Discussion of Theorem \ref{ratio}}
Let $\mathbb{T}= \left\{ t \in \mathbb{C} \given |t| = 1 \right\}$, and let $\sigma(t) = \sum_{i\in \mathbb{Z}} d_i t^i =\exp \left( \sum_{i\in \mathbb{Z}}c_i t^i \right)=e^{ f(t)}$ be a function on $\mathbb{T}$. 

We will always assume $f(t^{-1})=f(t)$ (i.e. $c_i=c_{-i}$). 

We will also need two extra conditions:
\begin{description}
\item[Condition (A)] $$\sum |c_i| < \infty$$
\item[Condition (B)] $$\sum |i| |c_i|^2 < \infty$$
\end{description}
Those conditions were already relevant to the work of Bump and Diaconis~\cite{BD}, and the whole field of Toeplitz matrices\footnote{
The book by B\"ottcher and Silbermann \cite{BS} gives a very clear introduction to the analytic theory of Toeplitz matrices. Theorem 5.2 in \cite{BS} uses those conditions. Sets of functions satisfying Conditions~(A) and (B) are denoted $W(\mathbb{T})$ and $B_2^{1/2}(\mathbb{T})$ respectively.}.

One can define a class function $\fatPhi_{n,f}(g)$ on $G$ as $$\fatPhi_{n,f}(g)=e^{nc_0}\exp\left(\sum_{i> 0} c_i \myPbdg{(i)}\right).$$
A possibly more intuitive definition (but only valid when $G=\symp{2n}$ or $G=\sorth{2n}$) is $\fatPhi_{n,f}(g)=\prod_{k=1}^n \sigma(t_k)$, where the product is taken over half of the eigenvalues of $g$, one in each conjugate pair. The symmetry condition $f(t^{-1})=f(t)$ guarantees that $\fatPhi_{n,f}$ is independent of the chosen subset of eigenvalues. When $G=\sorth{2n+1}$, the product expression becomes slightly more complicated because of the eigenvalue~1.

\begin{thm}
\label{ratio}
Assume that $f$ satisfies Condition~(A). 
For simplicity of
notation, take $ \chiGbd{\gamma}{}= \chisorthoddbd{\gamma}{}$(resp. $\chisympbd{\gamma}{}$, $\chisorthevenbd{\gamma}{}$) if $G=\sorth{2n+1}$
(resp. $\symp{2n},\sorth{2n}$). Then
$$
\lim_{n \rightarrow \infty} \frac{\esp_G \chiGbd{\gamma}{}\fatPhi_{n,f}}{\esp_G \fatPhi_{n,f}} = R(\gamma,(c_i)),
$$
with 
\begin{eqnarray*}
R(\gamma,(c_i)) &= &\sum_{\lambda \partition |\gamma|}  
\chi_\gamma(\lambda) \left(\prod_{i=1}^\infty 
\frac{c_i^{\lambda(i)}}{\lambda(i)!}\right) \\
&=&\biggl. s_\gamma \biggr\arrowvert_{p_i \coloneqq ic_i},
\end{eqnarray*}
where the last expression is a specialization for the Schur polynomial $s_\gamma$ when the value of the power polynomials is set using the Fourier coefficients $c_i$.
\end{thm}
We delay comments on this Theorem to page \pageref{remarks} and start with the proof.
\begin{proof}
As a first approximation to $\esp_G \chiGbd{\gamma} {}\fatPhi_{n,f}$, we 
will actually study $\esp_G \chiGbd{\gamma}{} \myPbd{\lambda}{}$ for $\lambda\partition k \le n$. It will be useful to split up $\lambda$ into subpartitions. To 
avoid confusion with notation previously used for partition parts 
($\lambda_1,\lambda_2,\cdots,\lambda_n$), we will use $\lambda_a \cup 
\lambda_b=\lambda$ in this proof only.

We start from the final equation in Proposition \ref{prime} and apply Frobenius Reciprocity to get
\begin{eqnarray*}
{\esp_G \chiGbd{\gamma}{} \myPbd{\lambda}{} }&=&  \left<
\chi_{\gamma}{} \otimes \restr^\sym{k-|\gamma|}_\centr{k-|\gamma|}
\sgn^\epsilon, \restr^\sym{k}_\SBgroup \myP{\lambda}{} \right>_\SBgroup
\\ 
&=&\frac
{z_{\lambda}}
{|\sym{|\gamma
|}||\centr{k-|\gamma|}|}\sum_{\substack{(\rho_a,\rho_b) \in
\SBgroup\\ \text{type}(\rho_a)=\lambda_a\partition
|\gamma|\\\text{type}(\rho_b) = \lambda_b \partition k-|\gamma|\\
\lambda_a \cup \lambda_b = \lambda}}
\chi_{\gamma}{(\rho_a)}\sgn^\epsilon(\rho_b),
\end{eqnarray*}
where  $\epsilon = 1$ when $G=\symp{2n}$ and
0 otherwise.
We now sum over conjugacy classes (i.e. cycle types) instead. The
correction factor for the $\rho_a$'s of type $\lambda_a$ will be
$\frac{|\sym{|\gamma|}|}{z_{\lambda_a}}= |\mathscr{C}_{\lambda_a}|$  , so
\begin{eqnarray*}
{\esp_G \chiGbd{\gamma}{} \myPbd{\lambda}{} }
&=&
\frac{z_\lambda}{|\centr{k-|\gamma|}|}
\sum_{\substack{\lambda_a\partition |\gamma| \\ \lambda_a \cup \lambda_b =
\lambda}}
\frac{\chi_{\gamma}{(\lambda_a)}\sgn(\lambda_b)^\epsilon}{z_{\lambda_a}}
|\centr{k-|\gamma|}\cap \mathscr{C}_{\lambda_b}|
.
\end{eqnarray*}
Observe from the proof of Theorem~\ref{traces}, with $\lambda$ replaced by
$\lambda_b$, that
$$
\esp_G \myPbd{{\lambda_b}}{} =
\frac{z_{\lambda_b}\sgn(\lambda_b)^\epsilon}{|\centr{k-|\gamma| }
|}|\centr{k-|\gamma| } \cap \mathscr{C}_{\lambda_b}|.
$$
The hypothesis $n\ge |\lambda_b|$ of Theorem~\ref{traces} is automatically satisfied since we already assume $n\ge |\lambda|$ and $\lambda = \lambda_a \cup \lambda_b$.

We now have the much simpler
$$
{\esp_G \chiGbd{\gamma}{} \myPbd{\lambda}{} }
 =  \sum_{\substack{\lambda_a\partition |\gamma| \\ \lambda_a \cup
\lambda_b = \lambda}}
\frac{z_\lambda}{z_{\lambda_a}z_{\lambda_b}}\,
\chi_{\gamma}{(\lambda_a)}\esp_G \myPbd{{\lambda_b}}
$$
or even
\begin{eqnarray}
\label{simpler}
{\esp_G \chiGbd{\gamma}{} \myPbd{\lambda}{} }
 =  \sum_{\substack{\lambda_a\partition |\gamma| \\ \lambda_a \cup
\lambda_b = \lambda}}
{  \frac{\lambda!}{\lambda_a!\lambda_b!}}
\chi_{\gamma}{(\lambda_a)}\esp_G \myPbd{{\lambda_b}}
\end{eqnarray}
where $\lambda!=\prod_{i\ge 1} (\lambda(i)!)$.

We can now deal with $\esp_G \chiGbd{\gamma}{}\fatPhi_{n,f}$. As in 
\emph{Toeplitz minors} \cite{BD}, absolute convergence is guaranteed by Condition~(A), 
the bound $|\tr(g^i)| \leq m$ when $g \in \unit{m}, \sorth{m}$ or 
$\symp{m}$ and compactness of those groups: 
% changed this to make sure the group symbols are as elsewhere. 
\begin{eqnarray*}
\esp_G\chiGbd{\gamma}{}\fatPhi_{n,f} &\leq& \int_G \max_{g \in 
G}(|\chiGbd{\gamma}{}|) \exp \left( \sum_{i \geq 0} |c_i| |\tr(g^i)| 
\right).
\end{eqnarray*}
We are thus allowed to permute sums and products in the full expansion of $\fatPhi_{n,f}$:
\begin{eqnarray*}
\esp_G \chiGbd{\gamma}{}\fatPhi_{n,f} &=&
e^{nc_0}\esp_G \chiGbd{\gamma}{} \exp \left( \sum_{i > 0} c_i \myPbd{(i)}{} \right)\\ 
&=&e^{nc_0}\esp_G \chiGbd{\gamma}{} \prod_{i=1}^\infty \sum_{j=0}^\infty \frac{(c_i \myPbd{(i)}{})^j}{j!}\\
&=&e^{nc_0}\esp_G \chiGbd{\gamma}{} \sum_{(\alpha_i)} \prod_{i=1}^\infty \frac{(c_i \myPbd{(i)}{})^{\alpha_i}}{\alpha_i!}\\
&=&e^{nc_0}\esp_G \chiGbd{\gamma}{} \sum_{(\alpha_i)} \prod_{i=1}^\infty \frac{c_i^{\alpha_i}}{\alpha_i!}\myPbd{(i^{\alpha_i})}{}\\
&=&e^{nc_0}\mathlarger{\sum_{\substack{(\alpha_i)\\\lambda \coloneqq (i^{\alpha_i})}}} \left(\prod_{i=1}^\infty \frac{c_i^{\alpha_i}}{\alpha_i!} \right) \esp_G \chiGbd{\gamma}{}  \myPbd{\lambda}{},
\end{eqnarray*}
From this definition of $\lambda$, we observe that $\lambda(j)=\alpha_j$, which explains the notation: $\alpha_j <> \lambda_j$ in general. 

Once $n\ge |\lambda|$, we are allowed to substitute for every term $\esp_G \chiGbd{\gamma}{}  \myPbd{\lambda}{}$ the r.h.s. of Equation~(\ref{simpler}). For a given $n$, this only applies for the terms at the head of the series, but any term in the series will eventually be substituted, when $n\ge |\lambda|$. Combined with absolute convergence, this guarantees the asymptotics
\begin{eqnarray*}
\esp_G \chiGbd{\gamma}{}\fatPhi_{n,f} &\overset{n \rightarrow \infty}{\sim}&e^{nc_0}{{\sum_{(\alpha_i)}}} \left( 
\left(\prod_{i=1}^\infty \frac{c_i^{\alpha_i}}{\alpha_i!} \right) \sum_{\substack{\lambda_a \partition |\gamma|\\ \lambda_a \cup \lambda_b = (i^{\alpha_i})=:\lambda}} \frac{\lambda!}{\lambda_a! \lambda_b!}\chi_\gamma(\lambda_a) \esp_G \myPbd{\lambda_b}{} \right).
\end{eqnarray*}
We now switch the sums, and change the index of one sum from $(\alpha_i)$ with $(i^{\alpha_i})=\lambda$ to $(\beta_i)$ with $(i^{\beta_i})=\lambda_b$. This implies $\lambda_a(j)+\beta_j=\lambda(j)=\alpha_j$. We get 
\begin{eqnarray*}
\esp_G \chiGbd{\gamma}{}\fatPhi_{n,f}&\overset{n \rightarrow \infty}{\sim}& e^{nc_0}{{\sum_{\lambda_a \partition |\gamma|} }} \left( \left(\frac{\chi_\gamma(\lambda_a)}{\lambda_a!}\prod_{i=1}^{\infty} c_i^{\lambda_a(i)} \right) {\sum_{(\beta_i)}}  \left(\prod_{i=1}^{\infty}\frac{ c_i^{\beta_i}}{\beta_i!}\right)\esp_G \myPbd{(i^{\beta_i})}{}\right)\\
&=&    R(\gamma,(c_i)) \esp_G \fatPhi_{n,f},
\end{eqnarray*}
and finally
\begin{eqnarray*}
\lim_{n \rightarrow \infty} \frac{\esp_G \chiGbd{\gamma}{}\fatPhi_{n,f}}{\esp_G \fatPhi_{n,f}} = R(\gamma,(c_i)) = \sum_{\lambda \partition |\gamma|} \chi_\gamma(\lambda) \left(\prod_{i=1}^\infty 
\frac{c_i^{\lambda(i)}}{\lambda(i)!}\right).
\end{eqnarray*}
The specialization expression now follows from the usual decomposition of power polynomials into Schur polynomials given by the character table of a symmetric group (see Sagan~\cite[Equation (4.23)]{Sagan}).
\end{proof}
\subsubsection*{Remarks. }
\label{remarks}
\begin{itemize}
\item As mentioned earlier, this ratio $R(\gamma,(c_i))$ already appears in Theorem 6 of Bump and Diaconis \cite{BD}, when $G=\unit{n}$. It is striking that this ratio is independent of the Cartan type of $G$. 
% ??? Link with universality in physics...
\item The authors went a bit further in \cite{BD} and modified the integrand using two characters (one of them appeared conjugated). There is no real need to do this here, as the characters $\chiGbd{\lambda}{}$ are real in the non-unitary cases, and we would just end up with a product of two characters. Koike and Terada \cite[Corollary 2.5.3]{KT} have shown that the multiplication rules are also essentially\footnote{This is only valid for $n\ge l(\mu) +l(\nu)$, and the case $G=\sorth{2n}$ is slightly different.} independent of the Cartan type of $G$, i.e. that 
$$
\chiGbd{\mu}{}\cdot\chiGbd{\nu}{} = \sum_\lambda c^\lambda_{\mu \nu} \chiGbd{\lambda}{}.
$$
This can be combined with Theorem~\ref{ratio} to show that there will also be an asymptotic ratio for $\frac{\esp_G \chiGbd{\mu}{}\chiGbd{\nu}{} \fatPhi_{n,f}}{\esp_G \fatPhi_{n,f}}$, independent of the Cartan type of $G$.
\item Johansson \cite[Theorem 3.8.i with $\eta = i$]{Jo} was the first to generalize the strong Szeg\"o limit theorem to all the classical groups. He found asymptotics for $\esp_G \fatPhi_{n,f}$ as $n\rightarrow \infty$. Bump and Diaconis~\cite{BD2} later found a new proof of Johansson's result that actually inspired our own work and an extension of this result. We state here a weaker version of Johansson's result in a style closer to our own. Note that this is the first time we need Condition~(B).
\begin{thm}[Johansson \cite{Jo}, Bump and Diaconis \cite{BD2}]
\label{Jo}
Let $f(t)=\sum_i c_i t^i$ satisfy Conditions~(A) and (B) in 
addition to the usual symmetry condition $f(t)=f(t^{-1})$. Then
\begin{eqnarray*}
\esp_{\sorth{2n+1}}\fatPhi_{n,f} & = & \exp\left( \sum_{i=1}^\infty \frac{ic_i^2}{2} - \sum_{i=1}^\infty c_{2i-1} + o(1)\right)\\
\esp_{\symp{2n}}\fatPhi_{n,f} & = & \exp\left( \sum_{i=1}^\infty \frac{ic_i^2}{2} - \sum_{i=1}^\infty c_{2i} + o(1)\right)\\
\esp_{\sorth{2n}}\fatPhi_{n,f} & = & \exp\left( \sum_{i=1}^\infty \frac{ic_i^2}{2} + \sum_{i=1}^\infty c_{2i}+ o(1) \right)\\
\end{eqnarray*}
\end{thm}

We can thus combine Theorems~\ref{ratio} and~\ref{Jo} to get the asymptotics for $\esp_G \chiGbd{\gamma}{}\fatPhi_{n,f}$, i.e. for the Haar measure twisted by a character of type $\chiGbd{\lambda}{}$.

\end{itemize}
\bibliographystyle{hplain}
\bibliography{references}
\end{document}